# Solving the Properties of Primes within Primes Sequence Based on the GRH is Tenable

**Kaida Shi**

Department of Mathematics, Zhejiang Ocean University,

Zhoushan City 316000, Zhejiang Province, China. (shikd@163.com)

**Abstract** By studying Dirichlet $L$ functional equation, the author obtained the relationship between Dirichlet $L$ functional equation and the prime property function, two new corollaries based on GRH is tenable and the equation of the prime property function. By using a geometric method, the author found a formula for solving the property $\chi(p)$ of primes within the prime sequence.

**2000MSC: 11M06, 11M20, 11M38.**

**Abbreviated title:    Solving Properties of Primes in Primes Sequence**

## 1. Introduction

After understanding there are infinite primes exist in the natural numbers, certainly, we must know that whether there are infinite primes exist in the arithmetic sequence:

$$l,\ l+q,\ l+2q,\ \cdots,\ l+dq,\ \cdots \tag{1}$$

(where, $q \geq 3,\ 1 \leq l < q,\ (l,\ q) = 1$)？L. Euler claimed that when $l = 1$, there are infinite primes exist in the arithmetic sequence (1). Also, A. M. Legendre considered certainly there are many primes exist in the arithmetic sequence (1). But, they didn't give their proofs. Although for special $l$ and $q$, people have proved many results in the elementary number theory. But, whether general conclusion is tenable is a very difficult conjecture. In 1837, Dirichlet proved for $q$ is a prime, this conjecture is true. After that, he use the **quadric form class number formula** (which was proved by himself) for general $q$, he found that this conjecture is also tenable. For fixing **a integer** whether belongs to the arithmetic sequence (1) ($d$ can be taken a negative value), he introduced a series of very important arithmetic function —— the property of the module $q$. Using it, we can select the **sub-sequence** from the govern **integer series**. When we investigate the **prime theorem** and other related number theory problems within the arithmetic sequence, the Dirichlet property is an essential and important tool[1].

I discover that for the properties of primes within the natural number sequence, people cannot solve them one by one. The aim of this paper is to give a formula for solving the properties of primes which are selected from the natural number sequence.

## 2. New corollaries based on GRH is tenable

We have Dirichlet $L$ functional equation:

$$\prod_{P}\left(1-\frac{\chi(p)}{p^s}\right)^{-1} = \sum_{n=1}^{\infty}\frac{\chi(n)}{n^s} = L(s,\chi),$$

It shows the relationship between all natural numbers (and they correspond the properties), and all primes (and they correspond the properties).

If take natural logarithm on two sides of the above identity, we obtain:

$$\log\prod_{P}\left(1-\frac{\chi(p)}{p^s}\right)^{-1} = \log L(s,\chi). \tag{2}$$

We can develop the left of the expression (2) as:

$$\log\prod_{P}\left(1-\frac{\chi(p)}{p^s}\right)^{-1} = -\sum_{P}\log\left(1-\frac{\chi(p)}{p^s}\right) =$$

$$= \sum_{P}\frac{\chi(p)}{p^s} + \frac{1}{2}\sum_{P}\left(\frac{\chi(p)}{p^s}\right)^2 + \frac{1}{3}\sum_{P}\left(\frac{\chi(p)}{p^s}\right)^3 + \cdots + \frac{1}{\mu}\sum_{P}\left(\frac{\chi(p)}{p^s}\right)^\mu + \cdots.$$

Considering the equation $L(s,\chi) = \sum_{n=1}^{\infty}\frac{\chi(n)}{n^s} = 0$, we can develop the right of the expression (2) as:

$$\log L(s,\chi) = (L(s,\chi)-1) - \frac{(L(s,\chi)-1)^2}{2} + \frac{(L(s,\chi)-1)^3}{3} - \cdots + (-1)^{\mu+1}\frac{(L(s,\chi)-1)^\mu}{\mu} - \cdots$$

$$= -(1-L(s,\chi)) - \frac{(1-L(s,\chi))^2}{2} - \frac{(1-L(s,\chi))^3}{3} - \cdots - \frac{(1-L(s,\chi))^\mu}{\mu} - \cdots.$$

As L. Euler used contrast method to proved audaciously $\sum_{n=1}^{\infty} \frac{1}{n^2} = \frac{\pi^2}{6}$ (**please refer to the Appendix below**), we contrast above two expressions, can obtain the following equation group:

$$\begin{cases} \omega_1(s,\chi) = \sum_p \frac{\chi(p)}{p^s} = -(1-L(s,\chi)); \\ \omega_2(s,\chi) = \sum_p \left(\frac{\chi(p)}{p^s}\right)^2 = -(1-L(s,\chi))^2; \\ \omega_3(s,\chi) = \sum_p \left(\frac{\chi(p)}{p^s}\right)^3 = -(1-L(s,\chi))^3; \\ \ldots\ldots\ldots\ldots\ldots\ldots\ldots\ldots\ldots; \\ \omega_\mu(s,\chi) = \sum_p \left(\frac{\chi(p)}{p^s}\right)^\mu = -(1-L(s,\chi))^\mu; \\ \ldots\ldots\ldots\ldots\ldots\ldots\ldots\ldots\ldots. \end{cases} \quad (3)$$

From first equation in the equation group (3), we can obtain:

$$\left(\sum_p \frac{\chi(p)}{p^s}\right)^\mu = (-1)^\mu (1-L(s,\chi))^\mu$$
$$= (-1)^\mu (-1)(-1)(1-L(s,\chi))^\mu$$
$$= (-1)^{\mu+1}(-(1-L(s,\chi))^\mu) = (-1)^{\mu+1} \sum_p \left(\frac{\chi(p)}{p^s}\right)^\mu.$$

Namely, the above are: when $\mu$ is odd number,

$$\left(\sum_p \frac{\chi(p)}{p^s}\right)^\mu = \sum_p \left(\frac{\chi(p)}{p^s}\right)^\mu,$$

when $\mu$ is even number,

$$\left(\sum_p \frac{\chi(p)}{p^s}\right)^\mu = -\sum_p \left(\frac{\chi(p)}{p^s}\right)^\mu$$

Substituting above results into the equation group (3), we obtain:

$$\begin{cases} \sum_p \dfrac{\chi(p)}{p^s} = L(s,\chi) - 1; \\ -\left(\sum_p \dfrac{\chi(p)}{p^s}\right)^2 = -(L(s,\chi) - 1)^2 \\ \left(\sum_p \dfrac{\chi(p)}{p^s}\right)^3 = (L(s,\chi) - 1)^3; \\ \cdots\cdots\cdots\cdots\cdots\cdots\cdots\cdots; \\ (-1)^{\mu+1}\left(\sum_p \dfrac{\chi(p)}{p^s}\right)^\mu = (-1)^{\mu+1}(L(s,\chi) - 1)^\mu; \\ \cdots\cdots\cdots\cdots\cdots\cdots\cdots\cdots. \end{cases}$$

This shows all the equations within above equation group are all identical.

So, we obtain:

$$\sum_p \frac{\chi(p)}{p^s} = L(s,\chi) - 1. \tag{4}$$

If $\omega(s,\chi) = \sum_p \dfrac{\chi(p)}{p^s}$ is called as the **prime property function**, then, the equation (4) shows the relationship between **Dirichlet $L$ function** $L(s,\chi) = \sum_{n=1}^{\infty} \dfrac{\chi(n)}{n^s}$ and the **prime property function**

$$\omega(s,\chi) = \sum_p \frac{\chi(p)}{p^s}.$$

**Corollary 1** If GRH is tenable, namely

$$L(s,\chi) = \sum_{n=1}^{\infty} \frac{\chi(n)}{n^s} = 0, \ (s = \frac{1}{2} + it, \ t \neq 0),$$

substituting the equation (4) into (3), we obtain:

$$\sum_p \frac{\chi(p)}{p^s} = -1.$$

**Corollary 2** Considering the **set of primes** add the **set of non-primes** equals to the set of natural numbers, we have:

$$\begin{cases} \omega(s,\chi) = \sum_p \dfrac{\chi(p)}{p^s} = -1; \\ L(s,\chi) = \sum_{n=1}^{\infty} \dfrac{\chi(n)}{n^s} = 0; \quad (s = \dfrac{1}{2} + it, \ t \neq 0) \\ \lambda(s,\chi) = \sum_c \dfrac{\chi(c)}{c^s} = 1; \end{cases}$$

Where, $\omega(s,\chi) = \sum_p \dfrac{\chi(p)}{p^s} = -1$ is called as the **equation of prime property function**, $\lambda(s,\chi) = \sum_c \dfrac{\chi(c)}{c^s} = 1$ is called as **equation of non-prime property function**.

## 3. Solving the properties of primes within the prime sequence based on GRH is tenable

We have known that there are infinite primes exist in the sequence of natural numbers. Denoting the property of prime $p$ as $\chi(p) = \alpha_p + i\beta_p$, we have:

$$\omega(s,\chi) = \sum_p \frac{\chi(p)}{p^{\frac{1}{2}+it}} = \sum_p \frac{\alpha_p + i\beta_p}{p^{\frac{1}{2}+it}} = \sum_p \frac{\alpha_p}{\sqrt{p}} e^{-it\log p} + i \cdot \sum_p \frac{\beta_p}{\sqrt{p}} e^{-it\log p}$$

$$= \sum_p \frac{\alpha_p}{\sqrt{p}} (\cos(t\log p) - i\sin(t\log p)) + i \cdot \sum_p \frac{\beta_p}{\sqrt{p}} (\cos(t\log p) - i\sin(t\log p))$$

$$= \sum_p \left( \frac{\alpha_p}{\sqrt{p}} \cos(t\log p) + \frac{\beta_p}{\sqrt{p}} \sin(t\log p) \right) + i \cdot \sum_p \left( \frac{\beta_p}{\sqrt{p}} \cos(t\log p) - \frac{\alpha_p}{\sqrt{p}} \sin(t\log p) \right) = -1.$$

Thus, we obtain

$$\omega(s,\chi) = \sum_p \frac{1}{\sqrt{p}} (\alpha_p \cos(t\log p) + \beta_p \sin(t\log p)) = -1.$$

Writing above expression as the form of the inner product between two infinite-dimensional vectors, we obtain:

$$\omega(s,\chi) = \sum_p \frac{1}{\sqrt{p}} (\alpha_p \cos(t \log p) + \beta_p \sin(t \log p)) = \left(\frac{1}{\sqrt{p_1}}, \frac{1}{\sqrt{p_2}}, \cdots, \frac{1}{\sqrt{p_m}}, \cdots\right) \bullet$$
$$\bullet (\alpha_{p_1} \cos(t_1 \log p_1) + \beta_{p_1} \sin(t_1 \log p_1), \ \alpha_{p_2} \cos(t_2 \log p_2) + \beta_{p_2} \sin(t_2 \log p_2), \ \cdots,$$
$$\alpha_{p_m} \cos(t_m \log p_m) + \beta_{p_m} \sin(t_m \log p_m), \ \cdots) = -1.$$

On the other hand, we can obtain:

$$\omega(s,\chi) = \left(\frac{1}{\sqrt{p_1}}, \frac{1}{\sqrt{p_2}}, \cdots, \frac{1}{\sqrt{p_m}}, \cdots\right) \cdot \left(\frac{-\sqrt{p_1}}{1 \cdot 2}, \frac{-\sqrt{p_2}}{2 \cdot 3}, \cdots, \frac{-\sqrt{p_m}}{m(m+1)} \cdots\right) = -1$$

Considering the **correspondence** between above expressions, we obtain:

$$\alpha_{p_m} \cos(t_m \log p_m) + \beta_{p_m} \sin(t_m \log p_m) = \frac{-\sqrt{p_m}}{m(m+1)}.$$

Because $|\chi(p_m)| = 1$, namely $\sqrt{\alpha_{p_m}^2 + \beta_{p_m}^2} = 1$, suppose that $\alpha_{p_m} = \sin u$ and $\beta_{p_m} = \cos u$, we obtain:

$$\sin u \cos(t_m \log p_m) + \cos u \sin(t_m \log p_m) = \frac{-\sqrt{p_m}}{m(m+1)},$$

i.e.

$$\sin(u + t_m \log p_m) = \frac{-\sqrt{p_m}}{m(m+1)}.$$

Therefore, we obtain:

$$\sin(\arcsin \alpha_{p_m} + t_m \log p_m) = \frac{-\sqrt{p_m}}{m(m+1)}$$

or

$$\sin(\arccos \beta_{p_m} + t_m \log p_m) = \frac{-\sqrt{p_m}}{m(m+1)}.$$

From

$$\arcsin \alpha_{p_m} = \arcsin\left(\frac{-\sqrt{p_m}}{m(m+1)}\right) - t_m \log p_m,$$

and

$$\arccos \beta_{p_m} = \arcsin\left(\frac{-\sqrt{p_m}}{m(m+1)}\right) - t_m \log p_m,$$

we obtain:

$$\begin{cases} \alpha_{p_m} = \sin\left(\arcsin\left(\frac{-\sqrt{p_m}}{m(m+1)}\right) - t_m \log p_m\right); \\ \beta_{p_m} = \cos\left(\arcsin\left(\frac{-\sqrt{p_m}}{m(m+1)}\right) - t_m \log p_m\right). \end{cases} \quad (5)$$

This is the formula for solving the real part $\alpha_{p_m}$ and imagine part $\beta_{p_m}$ of the property $\chi(p_m)$ of prime $p_m$ within the prime sequence.

About the formula for solving the imagine part $t_m$ (which corresponds to prime $p_m$) of non-trivial zero of the Riemann $\zeta$ function $\zeta(s)$, I have given it in my another paper entitled *New results based ob Riemann hypothesis is tenable*:

$$t_m = \frac{\arccos\left(\frac{-\sqrt{p_m}}{m(m+1)}\right)}{\log p_m},$$

So, (5) can be rewrite as:

$$\begin{cases} \alpha_{p_m} = \sin\left(\arccos\left(\frac{-\sqrt{p_m}}{m(m+1)}\right) - \arcsin\left(\frac{-\sqrt{p_m}}{m(m+1)}\right)\right); \\ \beta_{p_m} = \cos\left(\arccos\left(\frac{-\sqrt{p_m}}{m(m+1)}\right) - \arcsin\left(\frac{-\sqrt{p_m}}{m(m+)}\right)\right). \end{cases}$$

## 4. Conclusion

Although the permutation of primes has no law in the natural sequence[2~10], but every prime which is selected from the natural number sequence has its own subscript $m$, therefore, we can solve its property by using above formula.

Because the prime sequence is not arithmetic sequence, therefore, it has no certain relationship between the all the properties of primes. If we can obtain the permutation law of primes in the natural number sequence, then, probably, according to the formula $\chi(ab) = \chi(a)\chi(b)$, we can solve the properties of all **composite numbers** within the natural number sequence. But, unfortunately, the permutation law of primes within the natural number sequence is just a unsolved baffling problem. Hence, I hope the result of this paper will be of helpful to investigate deeply the property within the general arithmetic sequence.

# **Appendix**

# How was L. Euler proved $\sum_{n=1}^{\infty}\frac{1}{n^2}=\frac{\pi^2}{6}$ by using the contrast method ?

In the end of 17 century, Swiss mathematician J. Bernoulli wanted to know what the infinite series

$$1+\frac{1}{2^2}+\frac{1}{3^2}+\frac{1}{4^2}+\frac{1}{5^2}+\cdots\cdots$$

equals to ? If somebody solved the problem and announces him, he would appreciate him.

First, L. Euler (1707--1783) supposed that the equation

$$a_0 - a_1 x^2 + a_2 x^4 - \cdots + (-1)^n a_n x^{2n} = 0 \tag{1}$$

has $2n$ different roots:

$$\alpha_1,\ -\alpha_1,\ \alpha_2,\ -\alpha_2,\ \cdots,\ \alpha_n,\ -\alpha_n.$$

According to the factorization of the polynomial, he obtained:

$$a_0 - a_1 x^2 + a_2 x^4 - \cdots + (-1)^n a_n x^{2n}$$
$$= a_0\left(1-\frac{x^2}{\alpha_1^2}\right)\left(1-\frac{x^2}{\alpha_2^2}\right)\cdots\left(1-\frac{x^2}{\alpha_n^2}\right)$$
$$= a_0 - a_0\left(\frac{1}{\alpha_1^2}+\frac{1}{\alpha_2^2}+\cdots+\frac{1}{\alpha_n^2}\right)x^2 + \cdots.$$

Obviously,

$$a_1 = a_0\left(\frac{1}{\alpha_1^2}+\frac{1}{\alpha_2^2}+\cdots+\frac{1}{\alpha_n^2}\right). \tag{2}$$

Secondary, considering the equation:

$$\sin x = 0,$$

it equals to:

$$x - \frac{x^3}{3!} + \frac{x^5}{5!} - \frac{x^7}{7!} + \cdots = 0. \tag{3}$$

L.Euler considered the left hand of above equation has infinite terms and infinite-th power, therefore, it must has infinite roots:

$$0,\ \pi,\ -\pi,\ 2\pi,\ -2\pi,\ 3\pi,\ -3\pi,\ \cdots.$$

Removing the root 0 and use $x$ to divide two sides of (3), he obtained:

$$1 - \frac{x^2}{3!} + \frac{x^4}{5!} - \frac{x^6}{7!} + \cdots = 0 \tag{4}$$

The roots of above equation (4) are respectively:

$$\alpha_1 = \pi,\ \alpha_1' = -\pi,\ \alpha_2 = 2\pi,\ \alpha_2' = -2\pi,\ \alpha_3 = 3\pi,\ \alpha_3' = -3\pi,\ \cdots.$$

At the moment, L. Euler took (1) and (4) to **contrast audaciously**, and obtained:

$$a_0 - a_1 x^2 + a_2 x^4 - \cdots + (-1)^n a_n x^{2n}$$
$$= a_0 \left(1 - \frac{x^2}{\alpha_1^2}\right)\left(1 - \frac{x^2}{\alpha_2^2}\right)\cdots\left(1 - \frac{x^2}{\alpha_n^2}\right),$$
$$\frac{\sin x}{x} = 1 - \frac{x^2}{3!} + \frac{x^4}{5!} - \frac{x^6}{7!} + \cdots$$
$$= 1 \cdot \left(1 - \frac{x^2}{(\pm\pi)^2}\right)\left(1 - \frac{x^2}{(\pm 2\pi)^2}\right)\left(1 - \frac{x^2}{(\pm 3\pi)^2}\right)\cdots,$$

where

$$a_0 = 1, \ a_1 = \frac{1}{3!}, \ a_2 = \frac{1}{5!}, \ \ldots\ldots$$

Also, according to (2), he obtained:

$$\frac{1}{3!} = \frac{1}{\pi^2} + \frac{1}{4\pi^2} + \frac{1}{9\pi^2} + \cdots..$$

Taking $\pi^2$ to mulitplies two sides of above expression, he obtained:

$$\frac{\pi^2}{6} = 1 + \frac{1}{4} + \frac{1}{9} + \cdots,$$

namely

$$\frac{\pi^2}{6} = 1 + \frac{1}{2^2} + \frac{1}{3^2} + \cdots.$$

## The Properties $\chi(p_m)$ of Primes (first 180 primes, less than 1069)

| $m$ | $p_m$ | $\chi(p_m) = \alpha_{p_m} + \beta_{p_m} i$ | $m$ | $p_m$ | $\chi(p_m) = \alpha_{p_m} + \beta_{p_m} i$ |
|---|---|---|---|---|---|
| 1 | 2 | 0.0000000000-1.0000000000i | 45 | 197 | 0.9999080492-0.0135607190i |
| 2 | 3 | 0.8333333333-0.5527707984i | 46 | 199 | 0.9999148525-0.0130494317i |
| 3 | 5 | 0.9305555555-0.3661507313i | 47 | 211 | 0.9999170848-0.0128772499i |

| m | $p_m$ | $\chi(p_m) = \alpha_{p_m} + \beta_{p_m} i$ | m | $p_m$ | $\chi(p_m) = \alpha_{p_m} + \beta_{p_m} i$ |
|---|---|---|---|---|---|
| 4 | 7 | 0.9650000000-0.2622498808i | 48 | 223 | 0.9999193768-0.0126980302i |
| 5 | 11 | 0.9755555555-0.2197529477i | 49 | 227 | 0.9999243648-0.0122989667i |
| 6 | 13 | 0.9852607710-0.1710590926i | 50 | 229 | 0.9999295656-0.0118686113i |
| 7 | 17 | 0.9891581633-0.1468541046i | 51 | 233 | 0.9999337419-0.0115113760i |
| 8 | 19 | 0.9926697531-0.1208584350i | 52 | 239 | 0.9999370683-0.0112187094i |
| 9 | 23 | 0.9943209877-0.1064226175i | 53 | 241 | 0.9999411552-0.0108483202i |
| 10 | 29 | 0.9952066116-0.0977946843i | 54 | 251 | 0.9999430897-0.0106685212i |
| 11 | 31 | 0.9964416896-0.0842849878i | 55 | 257 | 0.9999458172-0.0104097418i |
| 12 | 37 | 0.9969592373-0.0779248296i | 56 | 263 | 0.9999483750-0.0101610685i |
| 13 | 41 | 0.9975244536-0.0703204417i | 57 | 269 | 0.9999507760-0.0099219707i |
| 14 | 43 | 0.9980498866-0.0624213410i | 58 | 271 | 0.9999537151-0.0096212109i |
| 15 | 47 | 0.9983680556-0.0571071418i | 59 | 277 | 0.9999557918-0.0094029000i |
| 16 | 53 | 0.9985672578-0.0535110425i | 60 | 281 | 0.9999580459-0.0091600431i |
| 17 | 59 | 0.9987398009-0.0501877478i | 61 | 283 | 0.9999604293-0.0088960536i |
| 18 | 61 | 0.9989569440-0.0456620630i | 62 | 293 | 0.9999615910-0.0087645051i |
| 19 | 67 | 0.9990720222-0.0430708083i | 63 | 307 | 0.9999622317-0.0086910962i |
| 20 | 71 | 0.9991950113-0.0401164469i | 64 | 311 | 0.9999640579-0.0084738155i |
| 21 | 73 | 0.9993159798-0.0369807036i | 65 | 313 | 0.9999659858-0.0082478581i |
| 22 | 79 | 0.9993828993-0.0351257829i | 66 | 317 | 0.9999675771-0.0080526230i |
| 23 | 83 | 0.9994552090-0.0330043212i | 67 | 331 | 0.9999681074-0.0079865059i |
| 24 | 89 | 0.9995055556-0.0314427163i | 68 | 337 | 0.9999693843-0.0078249869i |
| 25 | 97 | 0.9995408284-0.0303006989i | 69 | 347 | 0.9999702515-0.0077133737i |
| 26 | 101 | 0.9995901007-0.0286291900i | 70 | 349 | 0.9999717419-0.0075176699i |
| 27 | 103 | 0.9996395678-0.0268465000i | 71 | 353 | 0.9999729839-0.0073506132i |
| 28 | 107 | 0.9996754350-0.0254759632i | 72 | 359 | 0.9999740096-0.0072097307i |
| 29 | 109 | 0.9997119831-0.0239989764i | 73 | 367 | 0.9999748472-0.0070926041i |
| 30 | 113 | 0.9997386981-0.0228590352i | 74 | 373 | 0.9999757812-0.0069596724i |
| 31 | 127 | 0.9997418867-0.0227191549i | 75 | 379 | 0.9999766697-0.0068308101i |
| 32 | 131 | 0.9997650511-0.0216758538i | 76 | 383 | 0.9999776324-0.0066884061i |
| 33 | 137 | 0.9997823469-0.0208628569i | 77 | 389 | 0.9999784320-0.0065677578i |
| 34 | 139 | 0.9998036862-0.0198138613i | 78 | 397 | 0.9999790889-0.0064669765i |
| 35 | 149 | 0.9998122953-0.0193745758i | 79 | 401 | 0.9999799211-0.0063369886i |
| 36 | 151 | 0.9998297847-0.0184499741i | 80 | 409 | 0.9999805194-0.0062418673i |
| 37 | 157 | 0.9998411603-0.0178228542i | 81 | 419 | 0.9999810047-0.0061636231i |
| 38 | 163 | 0.9998515702-0.0172289764i | 82 | 421 | 0.9999818227-0.0060294420i |
| 39 | 167 | 0.9998627548-0.0165671853i | 83 | 431 | 0.9999822666-0.0059553779i |
| 40 | 173 | 0.9998713563-0.0160396627i | 84 | 433 | 0.9999830128-0.0058287293i |
| 41 | 179 | 0.9998792696-0.0155385416i | 85 | 439 | 0.9999835692-0.0057324872i |
| 42 | 181 | 0.9998890128-0.0148983950i | 86 | 443 | 0.9999841730-0.0056261621i |
| 43 | 191 | 0.9998932861-0.0146087803i | 87 | 449 | 0.9999846795-0.0055354133i |
| 44 | 193 | 0.9999015407-0.0140324262i | 88 | 457 | 0.9999850995-0.0054590089i |

| m | $p_m$ | $\chi(p_m) = \alpha_{p_m} + \beta_{p_m} i$ | m | $p_m$ | $\chi(p_m) = \alpha_{p_m} + \beta_{p_m} i$ |
|---|---|---|---|---|---|
| 89 | 461 | 0.9999856297-0.0053610071i | 135 | 761 | 0.9999954849-0.0030050324i |
| 90 | 463 | 0.9999861948-0.0052545447i | 136 | 769 | 0.9999955697-0.0029766873i |
| 91 | 467 | 0.9999866743-0.0051624727i | 137 | 773 | 0.9999956748-0.0029411666i |
| 92 | 479 | 0.9999869135-0.0051159413i | 138 | 787 | 0.9999957222-0.0029249807i |
| 93 | 487 | 0.9999872551-0.0050487317i | 139 | 797 | 0.9999957908-0.0029014552i |

| 94 | 491 | 0.9999876857-0.0049626991i | 140 | 809 | 0.9999958477-0.0028817524i |
|---|---|---|---|---|---|
| 95 | 499 | 0.9999880011-0.0048987370i | 141 | 811 | 0.9999959539-0.0028446741i |
| 96 | 503 | 0.9999883986-0.0048169236i | 142 | 821 | 0.9999960178-0.0028221283i |
| 97 | 509 | 0.9999887345-0.0047466789i | 143 | 823 | 0.9999961182-0.0027863198i |
| 98 | 521 | 0.9999889301-0.0047052899i | 144 | 827 | 0.9999962062-0.0027545575i |
| 99 | 523 | 0.9999893276-0.0046200267i | 145 | 829 | 0.9999963005-0.0027201071i |
| 100 | 541 | 0.9999893932-0.0046058109i | 146 | 839 | 0.9999963571-0.0026992331i |
| 101 | 547 | 0.9999896920-0.0045404719i | 147 | 853 | 0.9999963957-0.0026848812i |
| 102 | 557 | 0.9999899072-0.0044928208i | 148 | 857 | 0.9999964754-0.0026550459i |
| 103 | 563 | 0.9999901871-0.0044300901i | 149 | 859 | 0.9999965607-0.0026227004i |
| 104 | 569 | 0.9999904567-0.0043688029i | 150 | 863 | 0.9999966356-0.0025939812i |
| 105 | 571 | 0.9999907812-0.0042938996i | 151 | 877 | 0.9999966704-0.0025805299i |
| 106 | 577 | 0.9999910293-0.0042357204i | 152 | 881 | 0.9999967421-0.0025525989i |
| 107 | 587 | 0.9999912087-0.0041931515i | 153 | 883 | 0.9999968190-0.0025223064i |
| 108 | 593 | 0.9999914418-0.0041371966i | 154 | 887 | 0.9999968865-0.0024953935i |
| 109 | 599 | 0.9999916667-0.0040824730i | 155 | 907 | 0.9999968974-0.0024910188i |
| 110 | 601 | 0.9999919374-0.0040156023i | 156 | 911 | 0.9999969626-0.0024647031i |
| 111 | 607 | 0.9999921452-0.0039635331i | 157 | 919 | 0.9999970130-0.0024441659i |
| 112 | 613 | 0.9999923458-0.0039125773i | 158 | 929 | 0.9999970560-0.0024265169i |
| 113 | 617 | 0.9999925638-0.0038564568i | 159 | 937 | 0.9999971044-0.0024064806i |
| 114 | 619 | 0.9999927970-0.0037955249i | 160 | 941 | 0.9999971639-0.0023816538i |
| 115 | 631 | 0.9999929083-0.0037660673i | 161 | 947 | 0.9999972158-0.0023597380i |
| 116 | 641 | 0.9999930401-0.0037309068i | 162 | 953 | 0.9999972665-0.0023381562i |
| 117 | 643 | 0.9999932531-0.0036733886i | 163 | 967 | 0.9999972936-0.0023265452i |
| 118 | 647 | 0.9999934374-0.0036228675i | 164 | 971 | 0.9999973479-0.0023030933i |
| 119 | 653 | 0.9999935955-0.0035789669i | 165 | 977 | 0.9999973954-0.0022823643i |
| 120 | 659 | 0.9999937485-0.0035359442i | 166 | 983 | 0.9999974418-0.0022619444i |
| 121 | 661 | 0.9999939335-0.0034832518i | 167 | 991 | 0.9999974820-0.0022440927i |
| 122 | 673 | 0.9999940226-0.0034575776i | 168 | 997 | 0.9999975264-0.0022242383i |
| 123 | 677 | 0.9999941794-0.0034119048i | 169 | 1009 | 0.9999975552-0.0022112594i |
| 124 | 683 | 0.9999943143-0.0033721589i | 170 | 1013 | 0.9999976026-0.0021897242i |
| 125 | 691 | 0.9999944288-0.0033380117i | 171 | 1019 | 0.9999976441-0.0021706624i |
| 126 | 701 | 0.9999945248-0.0033091324i | 172 | 1021 | 0.9999976937-0.0021476725i |
| 127 | 709 | 0.9999946340-0.0032759618i | 173 | 1031 | 0.9999977244-0.0021333579i |
| 128 | 719 | 0.9999947258-0.0032478367i | 174 | 1033 | 0.9999977718-0.0021110213i |
| 129 | 727 | 0.9999948299-0.0032156115i | 175 | 1039 | 0.9999978095-0.0020930848i |
| 130 | 733 | 0.9999949452-0.0031795583i | 176 | 1049 | 0.9999978381-0.0020793690i |
| 131 | 739 | 0.9999950571-0.0031441731i | 177 | 1051 | 0.9999978824-0.0020579644i |
| 132 | 743 | 0.9999951787-0.0031052624i | 178 | 1061 | 0.9999979097-0.0020446286i |
| 133 | 751 | 0.9999952711-0.0030753391i | 179 | 1063 | 0.9999979521-0.0020238153i |
| 134 | 757 | -0.9999953735-0.0030418575i | 180 | 1069 | 0.9999979746-0.0020070932i |

### The Properties $\chi(c_\eta)$ of Composite Numbers (first 89, less than 121)

| $\eta$ | $c_\eta$ | $\chi(c_\eta) = a_{c_\eta} + b_{c_\eta} i$ | $\eta$ | $c_\eta$ | $\chi(c_\eta) = a_{c_\eta} + b_{c_\eta} i$ |
|---|---|---|---|---|---|
| 1 | 4 | -1.0000000000+0.0000000000i | 46 | 65 | 0.8542064722-0.5199339407i |
| 2 | 6 | -0.5527707984-0.8333333333i | 47 | 66 | -0.7223860797-0.6914899506i |
| 3 | 8 | 0.0000000000+1.0000000000i | 48 | 68 | -0.9891581633+0.1468541046i |
| 4 | 9 | 0.3888888888-0.9212846636i | 49 | 69 | 0.7697735079-0.6383171207i |
| 5 | 10 | -0.3661507313-0.9305555555i | 50 | 70 | -0.5973735392-0.8019631255i |
| 6 | 12 | -0.8333333333+0.5527707984i | 51 | 72 | 0.9212846636+0.3888888888i |
| 7 | 14 | -0.2622498808-0.9650000000i | 52 | 74 | -0.0779248296-0.9969592373i |
| 8 | 15 | 0.5730655308-0.8195095467i | 53 | 75 | 0.2332052936-0.9724275247i |
| 9 | 16 | 1.0000000000+0.0000000000i | 54 | 76 | -0.9926697531+0.1208584350i |
| 10 | 18 | -0.9212846636-0.3888888888i | 55 | 77 | 0.8837809268-0.4679009226i |
| 11 | 20 | -0.9305555555+0.3661507313i | 56 | 78 | -0.6871726268-0.7264941713i |
| 12 | 21 | 0.6592025906-0.7519653878i | 57 | 80 | 0.9305555555-0.3661507313i |
| 13 | 22 | -0.2197529477-0.9755555555i | 58 | 81 | -0.6975308636-0.7165547382i |
| 14 | 24 | 0.5527707984+0.8333333333i | 59 | 82 | -0.0703204417-0.9975244536i |
| 15 | 25 | 0.7318672839-0.6814471943i | 60 | 84 | -0.6592025906+0.7519653878i |
| 16 | 26 | -0.1710590926-0.9852607710i | 61 | 85 | 0.8666958863-0.4988368878i |
| 17 | 27 | -0.1851851859-0.9827036416i | 62 | 86 | -0.0624213410-0.9980498866i |
| 18 | 28 | -0.9650000000+0.2622498808i | 63 | 87 | 0.7752807973-0.6316167235i |
| 19 | 30 | -0.8195095467-0.5730655308i | 64 | 88 | 0.2197529477+0.9755555555i |
| 20 | 32 | 0.0000000000-1.0000000000i | 65 | 90 | -0.9996985133-0.024553662i |
| 21 | 33 | 0.6914899506-0.7223860794i | 66 | 91 | 0.9059164174-0.4234565442i |
| 22 | 34 | -0.1468541046-0.9891581633i | 67 | 92 | -0.9926697531+0.1064226175i |
| 23 | 35 | 0.8019631255-0.5973735392i | 68 | 93 | 0.7837777946-0.6210413581i |
| 24 | 36 | -0.3888888888+0.9212846636i | 69 | 94 | -0.0571071418-0.9983680556i |
| 25 | 38 | -0.1208584350-0.9926697531i | 70 | 95 | 0.8794819491-0.4759322441i |
| 26 | 39 | 0.7264941713-0.6871726268i | 71 | 96 | -0.5527707984-0.8333333333i |
| 27 | 40 | 0.3661507313+0.9305555555i | 72 | 98 | -0.5061422699-0.8624500000i |
| 28 | 42 | -0.7519653878-0.6592025906i | 73 | 99 | 0.1769276954-0.9842238514i |
| 29 | 44 | -0.9755555555+0.2197529477i | 74 | 100 | -0.7318672839+0.681447194i |
| 30 | 45 | 0.0245536624-0.9996985133i | 75 | 102 | -0.6691561682-0.7431218088i |
| 31 | 46 | -0.1064226175-0.9926697531i | 76 | 104 | 0.1710590926+0.9852607710i |
| 32 | 48 | 0.8333333333-0.5527707984i | 77 | 105 | 0.3380919562-0.9411130798i |
| 33 | 49 | 0.8624500000-0.5061422699i | 78 | 106 | -0.053511042-0.9985672578i |
| 34 | 50 | -0.6814471940-0.7318672839i | 79 | 108 | -0.1851851859+0.9827036416i |
| 35 | 51 | 0.7431218088-0.6691561682i | 80 | 110 | -0.5616927063i-0.8273459394i |
| 36 | 52 | -0.9852607710+0.1710590926i | 80 | 111 | 0.7877247941-0.6160273116i |
| 37 | 54 | -0.9827036416+0.1851851859i | 81 | 112 | 0.9650000000-0.2622498808i |
| 38 | 55 | 0.8273459394-0.5616927063i | 82 | 114 | -0.6494342145-0.7604177806i |
| 39 | 56 | 0.2622498808+0.9650000000i | 83 | 115 | 0.8863041999-0.4631035147i |
| 40 | 57 | 0.7604177806-0.6494342145i | 84 | 116 | -0.9952066116+0.0977946843i |
| 41 | 58 | -0.0977946843-0.9952066116i | 85 | 117 | 0.2255628478-0.9742286184i |
| 42 | 60 | -0.5730655308+0.8195095467i | 86 | 118 | -0.0501877478-0.9987398009i |
| 43 | 62 | -0.0842849878-0.9964416896i | 87 | 119 | 0.9160251562-0.4011208213i |
| 44 | 63 | 0.1336709845-0.9910257651i | 88 | 120 | 0.8195095467+0.5730655308i |
| 45 | 64 | -1.0000000000+0.0000000000i | 89 | 121 | 0.9034172839-0.4287624177i |